\documentclass[fleqn,11pt]{article}
\usepackage{amsmath,amssymb,amsthm,xcolor,framed,esint,dsfont}
\usepackage[english]{babel}
\usepackage[margin=3cm]{geometry}
\usepackage{csquotes}
\usepackage{hyperref}

\def\XXint#1#2#3{{\setbox0=\hbox{$#1{#2#3}{\int}$}
		\vcenter{\hbox{$#2#3$}}\kern-.5\wd0}}

\newcommand{\lt}{\left}
\newcommand{\rt}{\right}

\newcommand{\nn}{\nonumber}

\newcommand{\qd}{\quad}

\newcommand{\ep}{\epsilon}

\newcommand{\wt}{\widetilde}

\newcommand{\GI}{\mathcal{G}}

\newcommand{\HI}{\mathcal{H}}

\newcommand{\BI}{\mathcal{B}}

\newcommand{\R}{{\mathbb R}}

\def\Xint#1{\mathchoice
	{\XXint\displaystyle\textstyle{#1}}%
	{\XXint\textstyle\scriptstyle{#1}}%
	{\XXint\scriptstyle\scriptscriptstyle{#1}}%
	{\XXint\scriptscriptstyle\scriptscriptstyle{#1}}%
	\!\int}

\def\XXint#1#2#3{{\setbox0=\hbox{$#1{#2#3}{\int}$}
		\vcenter{\hbox{$#2#3$}}\kern-.5\wd0}}

\def\XXiint#1#2#3{{\setbox0=\hbox{$#1{#2#3}{\int}$}
		\vcenter{\hbox{$#2#3$}}\kern-.5\wd0}}

\newcommand{\e}{\epsilon}
\DeclareMathOperator{\supp}{supp}

\newtheorem{thm}{Theorem}[section]

\newtheorem{lem}[thm]{Lemma}

\theoremstyle{definition}

\newtheorem{rem}[thm]{Remark}

\title{On optimal regularity estimates for finite-entropy solutions of scalar conservation laws}

\date{}
\author{Xavier Lamy\footnote{Institut de Math\'ematiques de Toulouse, UMR 5219, Universit\'e de Toulouse, CNRS, UPS
		IMT, F-31062 Toulouse Cedex 9, France. Email: Xavier.Lamy@math.univ-toulouse.fr}
	\and Andrew Lorent\footnote{Department of Mathematical Sciences, University of Cincinnati, Cincinnati, OH 45221, USA. Email: lorentaw@uc.edu} 
	\and Guanying Peng\footnote{Department of Mathematical Sciences, Worcester Polytechnic Institute, Worcester, MA 01609, USA. Email: gpeng@wpi.edu}}

\begin{document}

\maketitle

\begin{abstract}
We consider finite-entropy solutions of  scalar conservation laws
$u_t +a(u)_x =0$, that is, bounded weak solutions whose entropy productions are locally finite Radon measures. 
Under the assumptions that the flux function $a$ is strictly convex (with possibly degenerate convexity) and $a''$ forms a doubling measure, we obtain a characterization of finite-entropy solutions in terms of an optimal regularity estimate involving a cost function first used by Golse and Perthame. 
\end{abstract}

\section{Introduction}

For any strictly convex $C^1$ flux function $a\colon\R\to\R$ we consider bounded weak solutions $u\colon [0,T]\times \R\to \R$ of the scalar conservation law 
\begin{equation}
\label{eq:scl}
u_t+\lt[a\lt(u\rt)\rt]_x=0\text{ in }\mathcal{D}'\lt( \lt (0,T\rt)\times\R\rt).
\end{equation}
It is well known that, on the one hand, smooth initial data evolving according to \eqref{eq:scl} may develop singularities in finite time, but on the other hand, there can be infinitely many weak solutions corresponding to a single initial datum. One way to restore well-posedness is the concept of entropy solution \cite{kruzkov70}.
For any convex or $C^2$ function $\eta\colon\R\to\R$, called entropy, and  entropy flux $q\colon\R\to \R$ such that $q'=\eta' a'$, the associated entropy production is the distribution
\begin{align*}
\mu_\eta =[\eta(u)]_t +[q(u)]_x.
\end{align*}
Entropy solutions are bounded weak solutions such that $\mu_\eta$ is a nonpositive measure for all convex entropies $\eta$, and for any bounded initial datum $u_0$ there exists a unique entropy solution defined for all positive times \cite{kruzkov70}.

Here we are interested in the larger class of solutions with finite entropy production:
\begin{align}\label{eq:finiteent}
\mu_\eta \in \mathcal M_{loc}([0,T]\times \R)\qquad\text{for all $C^2$ entropies }\eta.
\end{align}
This property does not ensure uniqueness of the initial value problem, but arises naturally in the study of some stochastic processes \cite{varadhan04,mariani10,BBMN10}, where large deviation principles are open for want of a better understanding of finite-entropy solutions (despite major recent progress in \cite{marconi-structure,marconi22rectifburgers}). 

In \cite[Proposition~2.3]{BBMN10} (see \cite[Appendix~B]{LPfacto} for a more detailed proof in a slightly different context) it is shown that \eqref{eq:finiteent} implies the existence of a locally finite Radon measure $m\in \mathcal M_{loc}([0,T]\times \R\times \R)$ such that
\begin{align}\label{eq:m}
\mu_\eta =\int \eta''(v) m(\cdot,\cdot,dv),
\end{align}
for all convex or $C^2$ entropy $\eta$. Then
 \cite[Theorem~4.1]{golseperthame13}
 implies
that  $u$ satisfies the regularity estimate
\begin{align}\label{eq:estGP}
& \sup_{|h|\leq\e}\frac{1}{|h|}\int_0^T \int_{-R}^{R} \chi(t,x)^2 \Delta(u(t,x),u(t,x+h))\, dx\, dt \nonumber \\
& \leq C(\chi) \left(1+|m|([0,T]\times [-R-\e,R+\e]\times [\inf u,\sup u])\right),
\end{align}
for any smooth cut-off function $\chi$   with support in $[0,T]\times [-R,R]$,  some constant $C(\chi)>0$, and the regularity cost $\Delta$ is given by
\begin{align}\label{eq:Delta}
\Delta(u_1,u_2)&=\frac 12 \int_{u_1}^{u_2}\int_{u_1}^{u_2} \lt|a'(v)-a'(w)\rt|\, dv\, dw \nonumber
\\
&=\int_{[u_1,u_2]}  (u_2-s)(s-u_1)  \, a''(ds).
\end{align}
The last equality is obtained by writing $|a'(v)-a'(w)|=\int_{[u_1,u_2]}(\mathbf 1_{v<s<w}+\mathbf 1_{w<s<v})\, a''(ds)$ and applying Fubini's theorem. 
Note that $[u_1,u_2]=[u_2,u_1]=\lbrace tu_1+(1-t)u_2\rbrace_{t\in [0,1]}$ and the integrand $(u_2-s)(s-u_1)$ is positive inside that segment, regardless of whether $u_1\leq u_2$ or $u_2\leq u_1$.

\begin{rem}
The explicit statement of \cite[Theorem~4.1]{golseperthame13} is actually a corollary of \eqref{eq:estGP}, but its proof does contain \eqref{eq:estGP}, which corresponds to (4.10) in the proof of \cite[Theorem~4.1]{golseperthame13}. The quantity $\Delta$ is defined in \cite[Lemma~4.3]{golseperthame13} by the formula
\begin{align*}
\Delta(u_1,u_2) & =\iint \mathbf 1_{v>w}(a'(v)-a'(w))(\mathcal M_{u_1}(v)-\mathcal M_{u_2}(v))(\mathcal M_{u_1}(w)-\mathcal M_{u_2}(w))\, dv\, dw,\\
\mathcal M_{u}(v)&=\mathbf 1_{0\leq v\leq u} -\mathbf 1_{u\leq v <0}.
\end{align*}
To see that this coincides with \eqref{eq:Delta}, first note that both expressions are  symmetric so it suffices to consider $u_1<u_2$. In the proof of \cite[Lemma~4.3]{golseperthame13} it is shown that
\begin{align*}
\Delta(u_1,u_2)= \int_{u_1}^{u_2}\int_w^{u_2} \lt|a'(v)-a'(w)\rt|\, dv \, dw
\end{align*}
which implies \eqref{eq:Delta}  by writing $a'(v)-a'(w)=\int_{[v,w]} a''(ds)$ and applying Fubini's theorem. 
\end{rem}

For instance, if $a(v)=|v|^{\beta+1}$ for some $\beta\geq 1$, then the regularity cost $\Delta$ admits the lower bound $\Delta(u_1,u_2)\gtrsim  |u_1-u_2|^{\beta+2}$. Hence in that case \eqref{eq:estGP} implies a local $B^{1/p}_{p,\infty}$ bound for $p=\beta+2$, in the $x$ direction, that is, $(t,x)\mapsto |u(t,x+h)-u(t,x)|/|h|^{1/p}$ is locally bounded in $L^p$, uniformly with respect to $h$.
In fact the same regularity is valid also in the $t$ direction \cite{golseperthame13}.
This local $B^{1/p}_{p,\infty}$ estimate is optimal in Besov regularity scales \cite{delelliswestdickenberg03}, but for $\beta >1$ it is
 strictly weaker than \eqref{eq:estGP}
 in regions where $u$ stays away from the degenerate value $u=0$. 
Loosely speaking, the regularity cost $\Delta$ takes into account that equation \eqref{eq:scl} regularizes more around values of $u$ where $a$ is more convex.
Therefore one can hope (as similar estimates in our recent work \cite{LLPon} for a generalized eikonal equation) that \eqref{eq:estGP} is optimal in the sense that a converse estimate is valid:
\begin{itemize}
\item If the left-hand side of \eqref{eq:estGP} is finite, does it imply that all entropy productions are finite \eqref{eq:finiteent} ?
\item Moreover,  are the entropy productions \eqref{eq:finiteent} controlled by the left-hand side of \eqref{eq:estGP} ? 
\end{itemize}

The second question can be answered rather easily if $a$ is $C^2$, thanks to the recent rectifiability result of \cite{marconi22rectifburgers}: under the \emph{a priori} knowledge that all entropy productions are finite, they are concentrated on a 1-rectifiable jump set and can be explicitly computed in terms of the traces of $u$ along that jump set. Elementary algebraic manipulation and a covering argument then provide the following estimate.

\begin{thm}\label{t:apriori}
	Assume that $a\in C^2(\R)$ is  strictly convex.  Let $u\in L^{\infty}([0,T]\times \R)$ be a weak solution of \eqref{eq:scl} such that $u$ has finite entropy production \eqref{eq:finiteent}. Then 
	for any open set $U\subset [0,T]\times\R$ we have the estimate
	\begin{align}\label{eq:mainest}
		|\mu_\eta|(U)\leq C_0\cdot   \sup_{I}|\eta''| \cdot \limsup_{\e\to 0} \sup_{|h|<\e}\frac{1}{|h|}\iint_U \Delta(u(t,x),u(t,x+h))\, dx\, dt,
	\end{align}
	for some absolute constant $C_0>0$, where $I=[\inf u,\sup u]$.
\end{thm}

Note that the \emph{a priori} estimate \eqref{eq:mainest} directly implies an estimate on $|m|(U\times\R)=|m|(U\times I)$ for the measure $m$ satisfying \eqref{eq:m}. In light of Theorem~\ref{t:apriori}, it is natural to reformulate the first question as follows:
 does finiteness of the right-hand side of \eqref{eq:mainest} imply finiteness of the left-hand side, that is, finite entropy production \eqref{eq:finiteent}?
 We provide a positive answer under a doubling assumption on the nonnegative measure $a''$.

\begin{thm}\label{t:main}
Assume that $a\in C^1(\R)$ is  strictly convex and that the nonnegative measure $a''$ is locally doubling, and let $u\in L^{\infty}([0,T]\times \R)$ be a weak solution of \eqref{eq:scl}. Assume that
\begin{align}\label{eq:regDelta}
\limsup_{\e\to 0} \sup_{|h|<\e}\frac{1}{|h|}\int_0^T \int_{-R}^R \Delta(u(t,x),u(t,x+h))\, dx\, dt <\infty,
\end{align}
for all $R>0$, then $u$ has finite entropy production \eqref{eq:finiteent}.
\end{thm}

Theorem~\ref{t:main} provides a full converse to the estimate \eqref{eq:estGP} proved in \cite{golseperthame13}, under the assumption that $a''$ is locally doubling (this is satisfied in particular if $a$ is analytic, see e.g. \cite[Lemma~25]{LLPon}). 
The proof of Theorem~\ref{t:main} also provides the estimate \eqref{eq:mainest} even when $a$ is not $C^2$, 
but with a constant depending on the doubling property of  $a''$.
More precisely, in the proof of Theorem~\ref{t:main} we obtain 
	\begin{align*}
		|\mu_\eta|(U)\leq C_0\cdot   \sup_{I}|\eta''| \cdot \limsup_{\e\to 0} \sup_{|h|<\e}\frac{1}{|h|}\iint_U \widehat\Delta(u(t,x),u(t,x+h))\, dx\, dt,
	\end{align*}
for some absolute constant $C_0$ and slightly different regularity cost $\widehat \Delta$ \eqref{eq:hatDelta}, and then check that $\widehat\Delta \leq C \Delta$ for some $C>0$ depending on the doubling constant of $a''$ on $I$.

Note that in the case $a(v)=|v|^{\beta+1}$ for some $\beta > 1$, 
the statement of Theorem~\ref{t:main} would not be valid with \eqref{eq:regDelta} replaced by a local $B^{1/p}_{p,\infty}$ bound for $p=\beta+2>3$. Indeed, for a solution taking values for instance in $[1,2]$ where $a$ is uniformly convex, $B^{1/3}_{3,\infty}$ regularity (in the $x$ direction) would be needed to ensure \eqref{eq:finiteent} (see the examples in \cite{delelliswestdickenberg03}).

It is also interesting to remark that,  if the limit \eqref{eq:regDelta} is zero, then all entropy productions vanish. In our particular context this provides a very precise regularity threshold for Onsager-type statements in the spirit of \cite{bardos19}, and a generalization of \cite[Theorem~2]{DeI} where $a(v)=v^2/2$ is considered.

The proof of Theorem~\ref{t:main} relies, as in \cite{DeI,bardos19}, on good estimates of the commutator $[a(u)]_\e-a(u_\e)$, where the subscript $\e$ denotes regularization at scale $\e$.
However, if the convexity of $a$ degenerates (e.g. $a(v)=|v|^{\beta+1}$ for some $\beta>1$), 
our regularity requirement \eqref{eq:regDelta} is strictly weaker than the local $B^{1/3}_{3,\infty}$ 
regularity that is needed in order to directly use (as done e.g. in \cite[Proposition~3.10]{GL}) the estimates of \cite[Theorem~2]{DeI}.
As noted in \cite{bardos19} these estimates are valid for any $C^2$ function $a$ and not related to its convexity. 
Here we  take instead full advantage of the convexity of $a$ in order to obtain  finer bounds in terms of the regularity cost $\Delta$.
We do this by adapting  ideas of \cite{LLPon}, where a result analogous to Theorem~\ref{t:main} has been established for a class of generalized eikonal equations with degenerate convexity.

We do not know whether Theorem~\ref{t:main} is valid without the requirement that the nonnegative measure $a''$ is doubling, even though the \emph{a priori} estimate of Theorem~\ref{t:apriori} suggests that this requirement is superfluous.
In the next two sections we give the proofs of Theorems~\ref{t:apriori} and \ref{t:main}, respectively.

\noindent\bf  Acknowledgments. \rm  X. L. received support from ANR project ANR-18-CE40-0023. A. L. gratefully acknowledges the support of the Simons foundation, collaboration grant \#426900.

\section{Proof of Theorem~\ref{t:apriori}}

Let $u$ be a bounded weak solution to \eqref{eq:scl} with finite entropy production \eqref{eq:finiteent}.  
The proof of \cite[Theorem~1]{marconi22rectifburgers}, where $a(v)=v^2/2$ is considered, actually uses only the facts that:
\begin{itemize}
\item  $u$ solves a kinetic formulation 
\cite[(3)]{marconi22rectifburgers}, which is a consequence of finite entropy production,  
\item the flux function $a$ is $C^2$ (to construct a Lagrangian representation \cite[Theorem~1.2]{marconi-structure}),
\item and $a'$ is an increasing function (see \cite[Proposition~6]{marconi22rectifburgers} and Step~2 of \cite[Theorem~10]{marconi22rectifburgers}). 
\end{itemize}
Hence it applies in our setting: there exists an $\mathcal{H}^1$-rectifiable set $J_u$ such that all entropy productions $\mu_\eta$ are absolutely continuous with respect to $\mathcal H^1_{\lfloor J_u}$. More precisely, $u$ has strong traces on both sides of $J_u$
and for any entropy $\eta$ we have
\begin{equation*}
\mu_{\eta} =\left( (\eta(u^+)-\eta(u^-))\nu_t + (q(u^+)-q(u^-))\nu_x \right) \, \mathcal{H}^1_{\lfloor J_u}, 
\end{equation*}
where $\nu=(\nu_t,\nu_x)$ is the unit normal to $J_u$ and $u^{\pm}$ are the traces. 
The equation \eqref{eq:scl} also provides the Rankine-Hugoniot condition
\begin{align*}
(u^+-u^-)\nu_t + (a(u^+)-a(u^-))\nu_x =0\qquad\text{a.e. on }J_u,
\end{align*}
so $\mu_\eta$ can be rewritten as
\begin{align}\label{eq:muetarectif}
\mu_\eta &=c_\eta(u^+,u^-)\, \nu_x \, \mathcal{H}^1_{\lfloor J_u}\\
c_\eta(u^+,u^-)&=
q(u^+)-q(u^-)-\frac{a(u^+)-a(u^-)}{u^+-u^-}(\eta(u^+)-\eta(u^-)).
\nonumber
\end{align}
The crucial fact here is that the entropy cost $c_\eta$ is controlled by $\Delta$.

\begin{lem}\label{l:Deltacontrolsceta}
For any $\eta\in C^2(\R)$ and $u^\pm\in\R$ we have
\begin{align*}
|c_\eta(u^+,u^-)|\leq\frac 12 \left( \sup_{[u^-,u^+]}|\eta''|\right)   \Delta(u^+,u^-).
\end{align*}
\end{lem}
\begin{proof}[Proof of Lemma~\ref{l:Deltacontrolsceta}]
Since both sides of the estimate are symmetric in $(u^+,u^-)$ we may assume $u^-<u^+$. Using $q'=\eta'a'$ and Fubini's theorem we have the identities
\begin{align}
c_\eta(u^+,u^-)&=\int_{u^-}^{u^+}\eta'(t)\left(a'(t)-\frac{a(u^+)-a(u^-)}{u^+-u^-}\right)\, dt \nonumber\\
&=\frac{1}{u^+-u^-}\int_{u^-}^{u^+}\eta'(t) \int_{u^-}^{u^+}(a'(t)-a'(s))\, ds\, dt
\nonumber\\
&=\frac{1}{u^+-u^-}\int_{[u^-,u^+]} w_\eta(\tau)\, a''(d\tau),
\label{eq:weta}
\end{align}
where
\begin{align*}
w_\eta(\tau)&=\int_{u^-}^{u^+}\eta'(t)\int_{u^-}^{u^+}\left(\mathbf 1_{s<\tau<t} -\mathbf 1_{t<\tau<s}\right)\, ds \, dt \\
&=\int_{u^-}^{u^+}\eta'(t) \left(\mathbf 1_{t>\tau}(\tau-u^-)-\mathbf 1_{t<\tau}(u^+-\tau)\right) \, dt.
\end{align*}
Since the second factor in the integrand has zero average on $[u^-,u^+]$ we deduce
\begin{align*}
|w_\eta(\tau)|&=\left| \int_{u^-}^{u^+}(\eta'(t)-\eta'(\tau)) \left(\mathbf 1_{t>\tau}(\tau-u^-)-\mathbf 1_{t<\tau}(u^+-\tau)\right) \, dt \right| \\
&\leq \left(\sup_{[u^-,u^+]}|\eta''|\right)
\int_{u^-}^{u^+}|t-\tau| \left(\mathbf 1_{t>\tau}(\tau-u^-)+\mathbf 1_{t<\tau}(u^+-\tau)\right) \, dt \\
& = \left(\sup_{[u^-,u^+]}|\eta''|\right)
\frac 12 (u^+-u^-)(\tau-u^-)(u^+-\tau).
\end{align*}
The last equality is obtained by directly calculating the integral. Plugging this into \eqref{eq:weta} we deduce
\begin{align*}
|c_\eta(u^+,u^-)|\leq  \frac 12 \left(\sup_{[u^-,u^+]}|\eta''|\right) \int_{[u^-,u^+]} (\tau-u^-)(u^+-\tau)\, a''(d\tau),
\end{align*}
and we recognize the definition \eqref{eq:Delta} of $\Delta(u^+,u^-)$ in the right-hand side.
\end{proof}

Theorem~\ref{t:apriori} follows from Lemma~\ref{l:Deltacontrolsceta} and the rectifiability of $J_u$ in \eqref{eq:muetarectif} by a covering argument similar to \cite[Lemma~32]{LLPon}.
We assume without loss of generality that $\sup_I|\eta''|\leq 1$ and $U\subset\subset[0,T]\times\R$. For general open $U\subset[0,T]\times\R$ we may approximate it by open sets $U_k\subset\subset[0,T]\times\R$.
Thanks to \eqref{eq:muetarectif} we have
\begin{align}\label{eq:muetaJ}
|\mu_\eta|(U)
&=\int_{J_u\cap U} |c_\eta(u^+,u^-)|\,|\nu_x| \, d\mathcal H^1 .
\end{align}
Further, for any $J'\subset J_u$ such that $\mathcal H^1(J')<\infty$ we have, on the one hand, thanks to Lemma~\ref{l:Deltacontrolsceta},
\begin{align}
\label{eqdcv1}
\int_{J'\cap U} |c_\eta(u^+,u^-)|\,|\nu_x| \, d\mathcal H^1 
\leq \frac 12 \int_{J'\cap U}\Delta(u^+,u^-)\, |\nu_x|\, d\mathcal H^1,
\end{align}
and, on the other hand, we will show
\begin{align}\label{eq:estDeltaJ}
\int_{J'\cap U}\Delta(u^+,u^-)\, |\nu_x|\, d\mathcal H^1
&\leq C_0\limsup_{\e\to 0} \sup_{|h|<\e}\frac{1}{|h|}
\iint_U \Delta(u(t,x),u(t,x+h))\, dx\, dt.
\end{align}
The proof of \eqref{eq:estDeltaJ} will follow as a consequence of the rectifiability of $J_u$ and the trace properties of $u$. Applying \eqref{eq:estDeltaJ}  to 
$J'_\delta=U\cap J_u\cap\lbrace  |c_\eta(u^+,u^-)\nu_x| >\delta\rbrace\subset J_u$and noting from \eqref{eq:muetaJ} that $\mathcal H^1(J'_\delta)\leq\delta^{-1}|\mu_\eta|(U)<\infty$, we deduce, thanks to \eqref{eqdcv1}, 
\begin{align*}
\int_{J'_\delta\cap U} |c_\eta(u^+,u^-)|\,|\nu_x| \, d\mathcal H^1 
\leq C_0  \limsup_{\e\to 0} \sup_{|h|<\e}\frac{1}{|h|}
\iint_U \Delta(u(t,x),u(t,x+h))\, dx\, dt.
\end{align*}
Letting $\delta\to 0$, the left-hand side converges to \eqref{eq:muetaJ}, and this proves the \emph{a priori} estimate \eqref{eq:mainest}.

To conclude the proof of Theorem~\ref{t:apriori} it remains to justify \eqref{eq:estDeltaJ}. The elementary building block is that \eqref{eq:estDeltaJ} is valid if $u$ is a pure jump: for constant values $u^\pm\in\R$ and a unit vector $\nu$, let $\psi^{u^\pm,\nu}\colon\R^2\to\R$ denote the pure jump from $u^-$ to $u^+$ across a line with unit normal $\nu$, namely
\begin{align*}
\psi^{u^\pm,\nu}(t,x)=u^- \mathbf 1_{(t,x)\cdot \nu<0} +u^+\mathbf 1_{(t,x)\cdot\nu >0},
\end{align*}
then  for any $r\geq |h|>0$ we claim
\begin{align}\label{eq:estDeltaJelem}
\int_{J_{\psi^{u\pm,\nu}}\cap B_r}\Delta(u^+,u^-) \,|\nu_x|\, d\mathcal H^1 \leq \frac{2}{|h|}\iint_{B_r}\Delta(\psi^{u^\pm,\nu}(t,x+h),\psi^{u^\pm,\nu}(t,x))\, dx\, dt.
\end{align}
To check \eqref{eq:estDeltaJelem}, simply use that the left-hand side is equal to $2r|\nu_x|\Delta(u^+,u^-)$, that $\Delta\geq 0$ and that the integrand in the right-hand side   is equal to $\Delta(u^+,u^-)$ in a region of two-dimensional measure $\geq r|h\nu_x|$ (the intersection of $B_r$ with a straight band of width $|h\nu_x|$).

We deduce \eqref{eq:estDeltaJ} from \eqref{eq:estDeltaJelem} via a covering argument similar to \cite[Lemma~32]{LLPon}, making use of the rectifiability of $J'$, the trace properties of $u$ and the Lipschitz quality of $\Delta$. We provide the details for the reader's convenience.

Let $\delta\in (0,1)$. There exists $\ep_0>0$ and a subset $\wt{J}\subset J'$ with $\mathcal{H}^1\lt(J'\cap U \setminus \wt{J}\rt)<\delta$ and $\widetilde J +B_{\e_0}\subset U$, such that for any $(t_0,x_0)\in \wt{J}$ and $0<r<\e_0$, denoting 
$u^\pm_0=u^\pm(t_0,x_0)$ and $\nu_0=\nu(t_0,x_0)$, we have
\begin{align}
	\label{eqmisbv11}
	\Xint{-}_{B_r(t_0,x_0)\cap J'} \left(\lt|u^{\pm}-u^{\pm}_0\rt| + |\nu-\nu_0|\right) d\mathcal{H}^1&<\delta,
	\nonumber\\
	\lt|\mathcal{H}^1\lt(B_r(t_0,x_0)\cap J'\rt)-2r \rt|&<\delta \, r,\nn\\
	\text{ and } \qd\frac{1}{\pi r^2}\iint_{B_r(t_0,x_0)} \lt|u-\psi^{u^\pm_0,\nu_0}_{t_0,x_0}\rt| dx\, dt &<\delta,
\end{align}
where  $\psi^{u^\pm_0,\nu_0}_{t_0,x_0}(t,x)=\psi^{u^\pm_0,\nu_0}(t-t_0,x-x_0)$ is the pure jump centered at $(t_0,x_0)$. 
Let $\ep\in (0,\ep_0/2)$. 
By Besicovitch's covering theorem \cite[Theorem~2.18]{ambrosio} there 
exists an absolute constant $Q\in \mathbb{N}$ and families $\BI_1, \BI_2,\dots, \BI_Q$ of pairwise disjoint balls in the set $\lt\{B_{\ep}(t,x)\colon (t,x)\in \wt{J}\rt\}$ such that 
\begin{equation*}
	\wt{J}\subset \bigcup_{k=1}^Q \bigcup_{B\in \BI_k} B.
\end{equation*}
We fix $k\in\lbrace 1,\ldots ,Q\rbrace$ and denote
 $\BI_{k}=\lbrace B_\e(t_j,x_j)\rbrace_{j=1,\ldots,p}$ for some $(t_j,x_j)\in\widetilde J$. We also write $u^\pm_j=u^\pm(t_j,x_j)$ and $\nu_j=\nu(t_j,x_j)$.

Note that $\Delta$ is Lipschitz on $I\times I$, with Lipschitz constant $L\lesssim |I| a''(I)$ thanks to its definition \eqref{eq:Delta}. Using the first two properties \eqref{eqmisbv11} of $\widetilde J$, we find
\begin{align*}
	\int_{J'\cap B_{\ep}(t_j,x_j)} \Delta\lt(u^{+},u^{-}\rt)\, |\nu_x| \, d\HI^1
	& \leq L\int_{J'\cap B_{\ep}(t_j,x_j)}\lt(|u^+-u_j^+|+|u^--u_j^-|\rt)|\nu_x|\,d\mathcal H^1\\
	&\qd+\Delta\lt(u_j^+,u_j^-\rt)\int_{J'\cap B_{\ep}(t_j,x_j)}|\nu_x-\nu_{j,x}|\,d\mathcal H^1\\
	&\qd+\Delta\lt(u_j^+,u_j^-\rt)|\nu_{j,x}|\mathcal H^1(J'\cap B_{\ep}(t_j,x_j))\\
	&\leq4\ep   \Delta\lt(u^{+}_j,u^{-}_j\rt) |\nu_{j,x}|+C\delta\e,
\end{align*}
for some constant $C=C(|I|, a''(I))$ depending on $|I|$ and $a''(I)$. Applying the elementary estimate \eqref{eq:estDeltaJelem} for pure jumps with $r=h=\e$, we deduce 
\begin{align*}
	&\int_{J'\cap B_{\ep}(t_j,x_j)} \Delta\lt(u^{+},u^{-}\rt) \,|\nu_x|\, d\HI^1
	\\
	& \leq 
	\frac{4}{\e}\iint_{B_\e(t_j,x_j)}
	\Delta(\psi_{t_j,x_j}^{u^\pm_j,\nu_j}(t,x+\e),\psi_{t_j,x_j}^{u_j^\pm,\nu_j}(t,x))\, dx\, dt
	+C\delta\e.
\end{align*}
And using the last property \eqref{eqmisbv11} of $\widetilde J$ we infer
\begin{align*}
	&\int_{J'\cap B_{\ep}(t_j,x_j)} \Delta\lt(u^{+},u^{-}\rt) \,|\nu_x|\, d\HI^1
	\\
	& \leq 
	\frac{4}{\e}\iint_{B_\e(t_j,x_j)}
	\Delta(u(t,x+\e),u(t,x))\, dx\, dt
	+(20\pi L+C)\delta\e.
\end{align*}
 Summing over $j=1,\ldots,p$ and over the families $\BI_1,\ldots,\BI_Q$ we obtain
\begin{align*}
	\int_{\wt{J}} \Delta\lt(u^{+},u^{-}\rt) \,|\nu_x |\, d\HI^1
	&\leq 4Q \frac{1}{\e} 
	\iint_U
	\Delta\lt(u(t,x+\e), u(t,x)\rt) \,dx\, dt +Q(20\pi L+C)\delta\, p\e.
\end{align*}
Noting from the properties \eqref{eqmisbv11} of $\widetilde J$ that 
\begin{align*}
	\mathcal{H}^1\lt(J'\cap U\rt)\geq \sum_{j=1}^p \mathcal{H}^1\lt(B_{\ep}(t_j,x_j)\cap J_u\rt)\geq p \ep,
\end{align*}
this implies
\begin{align*}
	\int_{\wt{J}} \Delta\lt(u^{+},u^{-}\rt) \,|\nu_x |\, d\HI^1
	&\leq 4Q \frac{1}{\e} 
	\iint_U
	\Delta\lt(u(t,x+\e), u(t,x)\rt) \,dx\, dt  \\
	&\quad 
	+Q(20\pi L+C)\delta\, \mathcal{H}^1\lt(J'\cap U\rt).
\end{align*}
Taking the limits $\e\to0$ and then $\delta\to 0$, we obtain \eqref{eq:estDeltaJ}.
\qed

\section{Proof of Theorem~\ref{t:main}}

We fix an entropy $\eta\in C^2(\R)$ and an entropy flux $q$ with $q'=\eta' a'$. 
The start of the proof is as in \cite[Theorem~2]{DeI}, we recall the argument for the reader's convenience.

We denote by a subscript $\e$ convolution at scale $\e$ in the $x$ variable:
\begin{align*}
u_\e(t,x)=\int u(t,z)\rho_\e(x-z)\, dz,
\end{align*}
where $\rho_\e(x)=\e^{-1}\rho(x/\e)$ for some smooth kernel $\rho\geq 0$ with $\supp\rho\subset [-1,1]$ and  $\int\rho=1$. We let 
\begin{align*}
\mu_\eta^\e =[\eta(u_\e)]_t +[q(u_\e)]_x,
\end{align*}
 and prove Theorem~\ref{t:main} by appropriately estimating $\mu_\eta^\e$. The regularized function $u_\e$ is pointwise differentiable with respect to $t$ and satisfies
\begin{align*}
u_{\e,t}=-[a(u)]_{\e,x},
\end{align*}
so we have
\begin{align*}
\mu_\eta^\e & =\eta'(u_\e)u_{\e,t} +q'(u_\e )u_{\e,x} \\
& = -\eta'(u_\e)[a(u)]_{\e,x} +\eta'(u_\e) [a(u_\e)]_x \\
& =\eta'(u_\e) \left[ a(u_\e)-[a(u)]_\e \right]_x \\
& = \left[\eta'(u_\e)  (a(u_\e)-[a(u)]_\e)\right]_x - \eta''(u_\e) u_{\e,x}(a(u_\e)-[a(u)]_\e).
\end{align*}
Testing this with a function $\psi\in C_c^\infty((0,T)\times \R)$ we obtain
\begin{align}\label{eq:muetaepspsi}
\left\langle \mu_\eta^\e,\psi\right\rangle &
 =-\iint \eta'(u_\e) (a(u_\e)-[a(u)]_\e) \psi_x\, dx\, dt\nonumber\\
&\quad
+ \iint \eta''(u_\e) u_{\e,x}([a(u)]_\e - a(u_\e)) \psi\, dx\, dt.
\end{align}
We have the convergences $u_\e\to u$ and $[a(u)]_\e\to a(u)$ a.e. and $u_\e$ is uniformly bounded, so by dominated convergence the left-hand side of \eqref{eq:muetaepspsi} converges to $\langle \mu_\eta,\psi\rangle$, and the first integral in the right-hand side of \eqref{eq:muetaepspsi} converges to 0. Hence we deduce
\begin{align}\label{eq:estmuetaeps}
\langle \mu_\eta,\psi\rangle 
\leq \|\psi\|_\infty  \sup_I |\eta''|\cdot  \limsup_{\e\to 0} \iint_{\supp\psi} |u_{\e,x}| ([a(u)]_\e - a(u_\e))\, dx\, dt.
\end{align}
Here recall that $I=[\inf u, \sup u]$, and note that 
\begin{align}\label{eq:jensen}
[a(u)]_\e - a(u_\e)\geq 0,
\end{align} 
by convexity of $a$ thanks to Jensen's inequality. Therefore it all boils down to estimating the right-hand side of \eqref{eq:estmuetaeps}, and this is where our proof needs to depart from \cite{DeI}.

We start by writing
\begin{align*}
&[a(u)]_{\e}(t,x) -a(u(t,x))=\int  \lt(a(u(t,z))-a(u(t,x))\rt)\rho_{\e}(x-z)\,dz\\
&\qquad=\int \lt(\int_{u(t,x)}^{u(t,z)}a'(\tau)\,d\tau\rt)\rho_{\e}(x-z)\,dz\\
&\qquad=\int \lt(\int_{u(t,x)}^{u(t,z)}(a'(\tau)-a'(u(t,x)))\,d\tau\rt)\rho_{\e}(x-z)\,dz\\
&\qquad\qd + a'(u(t,x))\int \lt(u(t,z)-u(t,x)\rt)\rho_{\e}(x-z)\,dz\\
&\qquad=\int \lt(\int_{u(t,x)}^{u(t,z)}(a'(\tau)-a'(u(t,x)))\,d\tau\rt)\rho_{\e}(x-z)\,dz + a'(u(t,x))\lt(u_\e(t,x)-u(t,x)\rt),
\end{align*}
hence
\begin{align}
\label{eqbdcv2}
	[a(u)]_{\e}(t,x) - a(u_{\e}(t,x)) 
	&=
[a(u)]_{\e}(t,x) -a(u(t,x))+a(u(t,x))- a(u_{\e}(t,x))	\nn\\	
	&= \int \lt(\int_{u(t,x)}^{u(t,z)}(a'(\tau)-a'(u(t,x)))\,d\tau\rt)\rho_{\e}(x-z)\,dz\nn\\
	&\qd+a(u(t,x))- a(u_{\e}(t,x)) + a'(u(t,x))\lt(u_\e(t,x)-u(t,x)\rt).
\end{align}
By convexity of $a$, we have
\begin{equation*}
	a(u(t,x))- a(u_{\e}(t,x)) + a'(u(t,x))\lt(u_\e(t,x)-u(t,x)\rt) \leq 0,
\end{equation*}
and applying this to \eqref{eqbdcv2} we deduce
\begin{equation}\label{eq:estcommut1}
 [a(u)]_{\e}(t,x) - a(u_{\e}(t,x)) \leq \int \lt(\int_{u(t,x)}^{u(t,z)}(a'(\tau)-a'(u(t,x)))\,d\tau\rt)\rho_{\e}(x-z)\,dz.
\end{equation}

To estimate this further, we define, for all $v\in\R$ and $r\geq 0$,
\begin{equation*}
	\GI_v(r)=\int_{v-r}^{v+r}\int_{v-r}^{v+r}|a'(\sigma)-a'(\tau)|\,d\sigma d\tau,
\end{equation*}
which satisfies
\begin{equation}\label{eq:G'}
	\GI_v'(r)=2\int_{v-r}^{v+r}\lt(|a'(v+r)-a'(\tau)|+|a'(v-r)-a'(\tau)|\rt)\,d\tau. 
\end{equation}
As $a'$ is strictly increasing, so is $\GI_v'$, and thus $\GI_v$ is strictly convex. Further, denoting by $g(x,z,t)=|u(t,z)-u(t,x)|$, we have
\begin{align*}
	\int_{u(t,x)}^{u(t,z)}|a'(\tau)-a'(u(t,x))|\,d\tau \leq \frac 12\GI_{u(t,x)}'(g(x,z,t)),
\end{align*}
and thus from \eqref{eq:estcommut1}, and recalling also \eqref{eq:jensen}, we infer
\begin{equation*}
	0\leq [a(u)]_{\e}(t,x) - a(u_{\e}((t,x)) \leq \frac 12\int_{B_{\e}}\GI_{u(t,x)}'(g(x,z,t))\rho_{\e}(x-z)\,dz.
\end{equation*}
Moreover we have
\begin{align*}
	\lt|u_{\ep,x}\lt(t,x\rt)\rt|\lesssim \ep^{-1} \Xint{-}_{[x-\e,x+\e]}  \lt|u(t,z)-u\lt(t,x\rt)\rt| dz = \ep^{-1} \Xint{-}_{[x-\e,x+\e]}  g(x,z,t)\, dz.
\end{align*}
Multiplying the last two estimates, we obtain
\begin{align}
	\label{eq:product}
	&\ep\lt([a(u)]_{\e}(t,x) - a(u_{\e}(t,x))\rt)\lt|u_{\ep,x}\lt(t,x\rt)\rt|\nn\\
	&\qquad \lesssim \Xint{-}_{[x-\e,x+\e]}  \Xint{-}_{[x-\e,x+\e]} \GI_{u(t,x)}'(g(x,z,t)) g(x,y,t)\,dy\,dz\nn\\
	&\qquad \leq  \Xint{-}_{[x-\e,x+\e]}\HI_{u(t,x)}\lt(\GI_{u(t,x)}'(g(x,z,t))\rt)\,dz+ \Xint{-}_{[x-\e,x+\e]}\GI_{u(t,x)}\lt(g(x,y,t)\rt)\,dy,
\end{align}
where $\HI_v(p)=\sup_{r\in\R} \lbrace pr-\GI_v(r)\rbrace$ is the Legendre transform of $\GI_v$. Using $\HI_v(p) = pr^*-\GI_v(r^*)$ where $r^*$ is characterized by $\GI_v'(r^*)=p$, we find that
\begin{align*}
\HI_{v}\lt(\GI_{v}'(r)\rt) &= r\GI_{v}'(r)-\GI_{v}(r)\\
&\leq r\GI_v'(r) \leq 8r^2 a''\lt([v-r, v+r]\rt).
\end{align*}
The last inequality follows from writing $a'(v+r)-a'(\tau)=a''([\tau,v+r])$ and $a'(\tau)-a'(v-r)=a''([v-r,\tau])$ in the explicit expression \eqref{eq:G'} of $\GI_v'$, and applying Fubini's theorem. Similarly we have
\begin{equation*}
	\GI_v(r)\leq 4r^2 a''\lt([v-r, v+r]\rt),
\end{equation*}
and plugging these bounds for $\HI_v(\GI_v'(r))$ and $\GI_v(r)$ into \eqref{eq:product}
 gives
\begin{align*}
	&\lt([a(u)]_{\e}(t,x) - a(u_{\e}(t,x))\rt)\lt|u_{\ep,x}\lt(t,x\rt)\rt|\nn\\
	& \lesssim \frac 1\e \Xint{-}_{[x-\e,x+\e]} g(x,z,t)^2 a''\lt([u(t,x)-g(x,z,t), u(t,x)+g(x,z,t)]\rt)\,dz, 
\end{align*}
where we recall that $g(x,z,t)=|u(t,z)-u(t,x)|$.
This implies
\begin{align*}
\lt([a(u)]_{\e}(t,x) - a(u_{\e}(t,x))\rt)\lt|u_{\ep,x}\lt(t,x\rt)\rt|
	& \lesssim\frac{1}{\e}  \Xint{-}_{[x-\e,x+\e]} \widehat\Delta(u(t,x),u(t,z))\, dz \\
	& = \frac{1}{\e}\Xint{-}_{[-\e,\e]} \widehat\Delta(u(t,x),u(t,x+h))\, dh,
\end{align*}
where
\begin{align}\label{eq:hatDelta}
\widehat \Delta (u_1,u_2) = |u_1-u_2|^2 a''([\min(u_1,u_2)-|u_1-u_2|,\max(u_1,u_2)+|u_1-u_2|]).
\end{align}
Integrating, we deduce
\begin{align*}
\iint_{\supp\psi}|u_{\e,x}| ([a(u)]_\e - a(u_\e))\, dx\, dt
& \lesssim\frac{1}{\e}\sup_{|h|<\e}\iint_{\supp\psi}\widehat\Delta(u(t,x),u(t,x+h))\, dx\, dt.
\end{align*}
Plugging this estimate into the bound \eqref{eq:estmuetaeps} for $\langle \mu_\eta,\psi\rangle$, we find
\begin{align*}
\langle \mu_\eta,\psi\rangle 
\lesssim \|\psi\|_\infty  \sup_I |\eta''|\cdot \limsup_{\e\to 0}\sup_{|h|<\e}\frac{1}{|h|} \iint_{\supp\psi} \widehat\Delta(u(t,x),u(t,x+h))\, dx\, dt.
\end{align*}
This is valid for any test function $\psi$ and implies in particular that $\mu_\eta$ is a locally finite Radon measure if the $\limsup$ in the right-hand side is finite.
It remains to show that, under the doubling assumption on $a''$, this $\limsup$ is controlled by \eqref{eq:regDelta}, thus concluding the proof of Theorem~\ref{t:main}.
 
Specifically, we claim
\begin{align}\label{eq:lowDelta}
\Delta(u_1,u_2)\geq C \widehat\Delta(u_1,u_2)\qquad\forall u_1,u_2\in I,
\end{align}
for some constant $C$ depending on the doubling constant of $a''$. 
To prove \eqref{eq:lowDelta} we may assume $u_1<u_2$. Letting $u_0=(u_1+u_2)/2$ and $r=|u_1-u_2|$, and recalling the explicit expression \eqref{eq:Delta} of $\Delta$, we have
\begin{align*}
\Delta(u_1,u_2)&=\int_{[u_1,u_2]}(s-u_1)(u_2-s) a''(ds) \geq \frac{r^2}{9}\int_{[u_0-r/6,u_0+r/6]} a''(ds).
\end{align*}
Thanks to the doubling property of $a''$ we deduce the lower bound
\begin{align*}
\Delta(u_1,u_2)\geq C r^2 a''([u_0-2r,u_0+2r]),
\end{align*}
which implies \eqref{eq:lowDelta} thanks to the explicit expression \eqref{eq:hatDelta} of $\widehat\Delta$, since $[u_0-2r,u_0+2r]$ contains $[\min(u_1,u_2)-r,\max(u_1,u_2)+r]$.
\qed

\bibliographystyle{acm}
\bibliography{aviles_giga}

\end{document}